\documentclass{amsart}
\usepackage{shortsalch}
\usepackage{xy}
\xyoption{all}
\title{Approximation of subcategories by abelian subcategories.}
\author{A. Salch}
\begin{document}
\maketitle

\section{Introduction.}

The earliest circulated versions of this paper date from 2010, and stem from ideas I had earlier, while I was a student. Several versions of this preprint have been made available, but the main idea in each of them is a very simple statement: {\em given an ideal $I$ of a commutative ring $R$ satisfying very mild hypotheses, the category of $L_0$-complete $R$-modules is the smallest abelian subcategory of $\Mod(R)$ containing all the $I$-adically complete $R$-modules and satisfying a few reasonable conditions}\footnote{You can take those conditions to be ``replete,'' ``exact,'' and ``full.'' It is not difficult to also prove variants of that result, replacing exactness with reflectivity, for example.}. As a slogan, ``the category of $L_0$-complete modules is the best abelian approximation to the category of $I$-adically complete modules.'' See Theorem \ref{main thm}, below, for a precise statement. 

In this short note, I prove the main theorem in close to its original level of generality. That level of generality is sufficient for every application I have ever had for the theorem, and it is also sufficient for every application I have ever seen anyone else have for the theorem. In this level of generality, the theorem is dramatically easier to prove than in the more general settings, as you can see from how short this note is! To me it seems that this short note is not worth sending to a journal, but I think it is worth having on the arXiv. None of the various versions of this paper were published, but in one version of this paper (the most general version) that I circulated some years ago, I believe there is a hypothesis missing from some of the statements, so I think it is worthwhile to post this short note, to serve as the simple, straightforward, and easily-seen-to-be-correct ``version of record'' for this theorem.

\section{The relevant definitions.}

\begin{definition}
Given an ideal $I$ in a commutative ring $R$, we write $\Lambda$ for the $I$-adic completion functor $\Lambda: \Mod(R)\rightarrow \Mod(R)$, i.e., $\Lambda(M) = \lim_n M/I^nM$. We write $L_0\Lambda$ for the zeroth left-derived functor $L_0\Lambda: \Mod(R)\rightarrow \Mod(R)$ of $\Lambda$.
\end{definition}
Since $I$-adic completion is, in general, not right exact, $L_0\Lambda(M)$ does not necessarily coincide with $\Lambda(M)$.

\begin{definition}
We say that an $R$-module $M$ is {\em $L_0$-complete} if the canonical map $\eta M: M \rightarrow L_0\Lambda M$ is an isomorphism. We write $L_0\Lambda\Mod(R)$ for the full subcategory of $\Mod(R)$ whose objects are the $L_0$-complete modules.
\end{definition}
See sections A.2 and A.3 of \cite{MR1601906} for an excellent introduction to $L_0$-completion and $L_0$-complete modules, including proofs of many basic properties.

\section{The theorem.}

Throughout, let $R$ be a commutative ring, and let $I$ be a {\em weakly pro-regular} ideal in $R$. Many early references on local homology and derived completion, such as \cite{MR1172439} and appendix A of \cite{MR1601906}, assumed that $R$ is Noetherian, or that $I$ is generated by a regular sequence. In \cite{MR1422312} and \cite{MR1973941}, it was established that weak pro-regularity of $I$ is sufficient for the proofs of most of the fundamental results in the area. Every weakly pro-regular ideal is finitely generated. 
Rather than reproduce the rather technical definition of weak pro-regularity here, I prefer to simply cite the result of \cite{MR1973941} which states that, in a Noetherian commutative ring, {\em every} ideal is weakly pro-regular. Consequently, in most practical situations, one knows that the ideals one encounters in examples are weakly pro-regular.

The papers \cite{MR3160712} and \cite{MR4352593} have valuable treatments of properties of $L_0\Lambda$ when $I$ is weakly pro-regular. For example, Theorem 3.9(a) of \cite{MR4352593} establishes that:
\begin{theorem}\label{old thm}
If $I$ is weakly pro-regular, then the full subcategory $L_0\Lambda\Mod(R)$ of $\Mod(R)$ is abelian, and the inclusion functor $\iota: L_0\Lambda\Mod(R)\rightarrow \Mod(R)$ is exact.
\end{theorem}
Theorem \ref{old thm} also appears in references that pre-date \cite{MR4352593}, although generally with stronger assumptions: for example, compare Theorem A.6 of \cite{MR1601906}, which is similar but includes the assumptions that $R$ is Noetherian and that $I$ is regular.

It is also straightforward (e.g. see Proposition 3.7 in \cite{MR4352593}) to prove that:
\begin{lemma}\label{basic lemma}
If $I$ is weakly pro-regular, then every $I$-adically complete $R$-module is $L_0$-complete.
\end{lemma}

With Theorem \ref{old thm} and Lemma \ref{basic lemma} in hand, we now have little trouble in proving the main theorem:
\begin{theorem}\label{main thm}
Let $I$ be a weakly pro-regular ideal in a commutative ring $R$. Then the category $L_0\Lambda\Mod(R)$ of $L_0$-complete modules is the unique smallest replete\footnote{Recall that a subcategory is {\em replete} if it is closed under isomorphisms, i.e., it contains every object isomorphic to one of its objects.} exact\footnote{Recall that a nonempty full abelian subcategory $\mathcal{A}$ of an abelian category $\mathcal{C}$ is {\em exact} if the inclusion functor $\mathcal{A}\hookrightarrow\mathcal{C}$ is exact.} full subcategory of $\Mod(R)$ containing all the $I$-adically complete $R$-modules.
\end{theorem}
\begin{proof}
Theorem \ref{old thm} establishes that $L_0\Lambda\Mod(R)$ is a replete exact full subcategory of $\Mod(R)$, while Lemma \ref{basic lemma} establishes that $L_0\Lambda\Mod(R)$ contains the $I$-adically complete $R$-modules. So suppose that $\mathcal{A}$ is a replete exact full subcategory of $\Mod(R)$ which contains the $I$-adically complete $R$-modules. We must then prove that $\mathcal{A}$ contains $L_0\Lambda\Mod(R)$. Suppose that $X$ is an $L_0$-complete $R$-module. Choose an exact sequence \[ P_1 \stackrel{d}{\longrightarrow} P_0\stackrel{}{\longrightarrow} X \rightarrow 0\] in $\Mod(R)$, with $P_0,P_1$ projective $R$-modules. Applying $\Lambda$, we have a short exact sequence
\begin{equation}\label{ses 1} 0 \rightarrow \im \Lambda d \rightarrow \Lambda P_0 \rightarrow L_0\Lambda X\rightarrow 0.\end{equation}
We know that $\Lambda P_0$ is in $\mathcal{A}$, since $\mathcal{A}$ contains all the $I$-adically complete $R$-modules.
If we can show that $\im \Lambda d$ is $I$-adically complete, then the short exact sequence \eqref{ses 1} exhibits $X\cong L_0\Lambda X$ as the quotient of an $I$-adically complete $R$-module by an $I$-adically complete submodule, i.e., as the quotient of an object of $\mathcal{A}$ by a subobject also in $\mathcal{A}$. Exactness and repleteness of $\mathcal{A}$ then gives us that $L_0\Lambda X$ is in $\mathcal{A}$ as well.

So all that is left is to show that the image of an $R$-module map with $I$-adically complete domain and $I$-adically complete codomain is also $I$-adically complete. Suppose $f: Y \rightarrow Z$ is a morphism of $R$-modules, with $Y,Z$ each $I$-adically complete. We have the commutative diagram
\[\xymatrix{
 Y \ar[r]^{\pi} \ar[d]^{\eta Y}_{\cong} 
  & \im f \ar[r]^i \ar[d]^{\eta \im f}
  & Z \ar[d]^{\eta Z}_{\cong} \\
 \Lambda Y \ar[r]^{\Lambda\pi} 
  & \Lambda \im f \ar[r]^{\Lambda i}
  & \Lambda Z 
}\]
where $\pi$ and $i$ are the canonical projection to the image and inclusion of the image, and $\eta: \id_{\Mod(R)}\rightarrow \Lambda$ is the canonical natural transformation sending each module to its $I$-adic completion.
Since $\Lambda$ preserves epimorphisms, $\Lambda \pi$ is epic, so $(\Lambda \pi)\circ \eta Y = (\eta(\im f))\circ \pi $ is epic, so $\eta(\im f)$ is epic. Meanwhile, $(\eta Z)\circ i = (\Lambda i) \circ \eta(\im f)$ is a composite of monomorphisms, so $\eta(\im f)$ is monic. So $\eta (\im f)$ is an isomorphism, i.e., $\im f$ is $I$-adically complete, as desired.
\end{proof}

\bibliography{/home/asalch/texmf/tex/salch}{}

\def\cprime{$'$} \def\cprime{$'$} \def\cprime{$'$} \def\cprime{$'$}
\begin{thebibliography}{1}

\bibitem{MR1422312}
Leovigildo Alonso~Tarr\'{\i}o, Ana Jerem\'{\i}as~L\'{o}pez, and Joseph Lipman.
\newblock Local homology and cohomology on schemes.
\newblock {\em Ann. Sci. \'{E}cole Norm. Sup. (4)}, 30(1):1--39, 1997.

\bibitem{MR1172439}
J.~P.~C. Greenlees and J.~P. May.
\newblock Derived functors of {$I$}-adic completion and local homology.
\newblock {\em J. Algebra}, 149(2):438--453, 1992.

\bibitem{MR1601906}
Mark Hovey and Neil~P. Strickland.
\newblock Morava {$K$}-theories and localisation.
\newblock {\em Mem. Amer. Math. Soc.}, 139(666):viii+100, 1999.

\bibitem{MR4352593}
Luca Pol and Jordan Williamson.
\newblock The homotopy theory of complete modules.
\newblock {\em J. Algebra}, 594:74--100, 2022.

\bibitem{MR3160712}
Marco Porta, Liran Shaul, and Amnon Yekutieli.
\newblock On the homology of completion and torsion.
\newblock {\em Algebr. Represent. Theory}, 17(1):31--67, 2014.

\bibitem{MR1973941}
Peter Schenzel.
\newblock Proregular sequences, local cohomology, and completion.
\newblock {\em Math. Scand.}, 92(2):161--180, 2003.

\end{thebibliography}
\bibliographystyle{plain}
\end{document}